# Characterization of count data distributions involving additivity and binomial subsampling

PEDRO PUIG[1] and JORDI VALERO[2]

[1]*Departament de Matemàtiques, Universitat Autònoma de Barcelona, 08193 Bellaterra, Spain. E-mail:* ppuig@mat.uab.es
[2]*Escola Superior d'Agricultura de Barcelona, Universitat Politècnica de Catalunya and Departament de Mathemàtiques, Universitat Autònoma de Barcelona. E-mail:* jordi.valero@upc.edu

In this paper we characterize all the $r$-parameter families of count distributions (satisfying mild conditions) that are closed under addition and under binomial subsampling. Surprisingly, few families satisfy both properties and the resulting models consist of the $r$th-order univariate Hermite distributions. Among these, we find the Poisson $(r = 1)$ and the ordinary Hermite distributions $(r = 2)$.

## 1. Introduction

An important problem in data analysis is how to choose an adequate family of distributions or statistical model to describe the values observed in any study. For this purpose the characterization theorems can be useful because, under general reasonable suppositions related to the nature of the experiment or the data, they allow us to reduce the possible set of distributions that can be used. One of these reasonable assumptions is that the model is closed under addition, that is, closed under convolutions. This means that the distribution of the sum of independent random variables with distributions belonging to the model also belongs to the same model. The property of closure under addition has widely been used to characterize families of distributions, the papers of Teicher [16] and Godambe and Patil [4] being particularly noteworthy. Other properties that can be utilized in characterization theorems jointly with additivity involve assumptions about the maximum likelihood estimators of the parameters. Puig [13] and Puig and Valero [14] characterize count families of distributions using several notions of additivity and assuming that the maximum likelihood estimator of the population mean is the sample mean.

In this paper we characterize count models using the additional notion of closure under binomial subsampling or under independent $p$-thinning. This property, studied in Wiuf and Stumpf [18], is properly defined in Section 2. The notion of binomial subsampling is quite natural in practice. Examples from mathematical biology can be found in Wiuf and Stumpf [18]. Roughly speaking, the property of closure under binomial subsampling means that if a count experiment is governed by a model with this property, the same family of distributions can be used to fit the data sets independently of the degree of loss of information due to the subsampling action.

From a result of Wiuf and Stumpf [18] one can directly establish the structure of the probability generating function of any count statistical model, with finite moments of all orders, that is closed under binomial subsampling. Corollary 1 in Section 2 extends this result to the first $r$ finite moments. Proposition 2 in Section 2 establishes that these models are constituted by mixed Poisson distributions and distributions of independent $p$-thinnings of random variables that do not come from any particular variable by independent $p$-thinning. Following Thedéen [17], these random variables that cannot be generated from an independent $p$-thinning of any random variable will be called *top inverses*. In order to illustrate all these results, three examples are given in Section 2.

The main result of this paper is stated in Theorem 1 in Section 3: all the $r$-parameter families of count distributions (satisfying mild conditions) that are closed under addition and under binomial subsampling are found. These are the $r$th-order univariate Hermite distributions of Milne and Westcott [11]. The case $r = 1$ corresponds to the Poisson model, and $r = 2$ is the ordinary Hermite family of distributions.

## 2. Count data models closed under binomial subsampling

We begin with a definition that can be found in Grandell [5], page 25.

**Definition 1.** *Let $X$ be a count random variable, let $\xi_1, \xi_2, \ldots$ be independent and identically distributed Bernoulli variables with probability of success $p$, for some $p \in (0, 1]$, and let $X$ be independent of the $\xi_i$. The count random variable*

$$X_p = \xi_1 + \xi_2 + \cdots + \xi_X \qquad (X_p = 0 \text{ if } X = 0)$$

*is called an independent $p$-thinning of $X$.*

The concept of $p$-thinning is widely used in point process theory, and it is also useful to describe count random variables coming from practical scenarios where a certain loss of information occurs. The random variable $X_p$ can be understood as a binomial subsampling of $X$. For instance, suppose that $X$ is the random variable which counts the number of seeds falling in a plot. Assume that each seed has the same probability $p$ of germinating and becoming a plant, and that for all seeds these actions are mutually independent and independent of the total number of seeds. Then the number of plants in the plot is $X_p$, a binomial subsampling of the number of seeds. The concept of *p*-thinning can be naturally generalized for multidimensional count random variables by a componentwise action. For instance, given a bivariate count random variable $\mathbf{X} = (Y, Z)$, an independent **p**-thinning ($\mathbf{p} = (r, s) \in (0, 1]^2$) can be defined as $\mathbf{X_p} = (Y_r, Z_s)$, where $Y_r$ and $Z_s$ are independent *r*- and *s*-thinnings. It is interesting to remark that if $\rho$ is the correlation coefficient of the two components of $\mathbf{X}$ and $\rho_\mathbf{p}$ the correlation coefficient of the components of its **p**-thinning, then $\rho_\mathbf{p}^2 \leq \rho^2$, and equality holds in the trivial case when $r = 1$, $s = 1$. Consequently, the **p**-thinning action diminishes the linear association of the components.

Let $\Phi_X(t)$ be the probability generating function (pgf) of a one-dimensional count random variable $X$. Recall that any pgf is a real analytic function, at least in $t \in (-1, 1)$. Moreover, if $X$ has finite moments of all orders, then $\Phi_X(t)$ is $C^\infty$ at $t = 1$ and $\Phi_X^{(k)}(1) = \mathrm{E}(X(X - 1) \cdots (X - k + 1)) = \mu_{(k)}$, that is, its *k*th factorial moment. When $X$ has finite moments of all orders, the *k*th derivative of $\log(\Phi_X(t))$ at $t = 1$ is known as the *k*th factorial cumulant of $X$, namely $\kappa_{(k)}$, and they are also quantities of interest for our purposes.

The pgf of an independent *p*-thinning or binomial subsampling of $X$ is (Grandell [5], page 25):

$$\Phi_{X_p}(t) = \Phi_X(1 - p(1 - t)). \tag{1}$$

The following lemma shows some elementary properties of an independent *p*-thinning that will be useful later on. The lemma follows immediately from (1).

**Lemma 1.** *Let $X$ be a count random variable with finite moments of all orders, and let $\mu_{(k)}$ and $\kappa_{(k)}$ be its $k$th factorial moment and cumulant, respectively. Its population mean $\mu_{(1)} = \kappa_{(1)}$ will be denoted by $\mu$. Let $\mu_{(k)}^*$ or $\kappa_{(k)}^*$ denote the $k$th factorial moment or cumulant for $X_p$, a binomial subsampling of $X$. Then $\mu_{(k)}^* = p^k \mu_{(k)}$ and $\kappa_{(k)}^* = p^k \kappa_{(k)}$. Moreover, the quantities $\mu_{(k)}/\mu^k$ and $\kappa_{(k)}/\mu^k$ are invariant under binomial subsampling, that is, $\mu_{(k)}/\mu^k = \mu_{(k)}^*/(\mu^*)^k$ and $\kappa_{(k)}/\mu^k = \kappa_{(k)}^*/(\mu^*)^k$ for any binomial subsampling.*

Using the notation $\eta_i = \kappa_{(i+1)}/\mu^{i+1}$ and the well-known relationships between the factorial cumulants and the ordinary cumulants (for instance, see Appendix 7 of Douglas [1]) the expressions below follow immediately:

$$\eta_1 = \frac{\sigma^2 - \mu}{\mu^2}, \qquad \eta_2 = \frac{\kappa_3 - 3\sigma^2 + 2\mu}{\mu^3}, \qquad \eta_3 = \frac{\kappa_4 - 6\kappa_3 + 11\sigma^2 - 6\mu}{\mu^4}, \tag{2}$$

where $\sigma^2$, $\kappa_3$ and $\kappa_4$ denote the variance and the third and fourth cumulants, respectively. Using the same technique we could find other quantities which are invariant under binomial subsampling involving higher-order cumulants. However, the cumulants of order more than 4 are not very useful from the statistical point of view.

Consider now a count statistical model, that is, a set of probability functions $\mathcal{F}$ of count random variables indexed by the parameter $\theta \in \Theta \subseteq \mathbb{R}^k$. Suppose that the model satisfies some minimum regularity conditions with respect to the parameters, that is, the model can be identified by a continuous function from $(0,1) \times \Theta \subset \mathbb{R}^{k+1}$ into $\mathbb{R}$, namely $\Phi(t; \theta)$, where, for each fixed value of $\theta_0 \in \Theta$, $\Phi(t; \theta_0)$ is a pgf. Moreover, we also assume that the parameter domain $\Theta$ is the widest possible. That is, if $\Phi^*(t; \theta)$ is any continuous function from $(0,1) \times \Theta^* \subset \mathbb{R}^{k+1}$ into $\mathbb{R}$, $\Theta \subset \Theta^*$, such that $\Phi^*(t; \theta) = \Phi(t; \theta)$ for all $\theta \in \Theta$, then for each $\theta_0 \in \Theta^* \setminus \Theta$, $\Phi^*(t; \theta_0)$ is not a pgf. The following definition, given in Wiuf and Stumpf [18], will be essential for our purposes.

**Definition 2.** *Let $\mathcal{F}$ be a count statistical model. It will be called closed under binomial subsampling or under independent $p$-thinning if, for any random variable $X$ with distribution belonging to the model, all its independent $p$-thinnings, for any $p \in (0,1]$, have distributions that also belong to the model.*

The Poisson distribution, the most frequently used statistical model to analyse count data, is closed under binomial subsampling. By way of explanation, for any Poisson random variable $X$, any independent $p$-thinning is also a Poisson random variable. However, the zero truncated Poisson distribution is not closed under independent $p$-thinning for any $p \in (0,1]$ except for the trivial case $p = 1$. Following Thedéen [17], a random variable $Y$ that is not an independent $p$-thinning of any random variable is called a *top inverse*. It is obvious that $Y$ is a top inverse if its distribution has at least one zero-gap, that is, if there exists $0 \leq k < n$ such that $P(Y = k) = 0$ while $P(Y = n) > 0$.

Models closed under binomial subsampling are very significant in practice. If a count experiment is ruled by one of these models, it means that the same family of distributions can be used to fit the data sets independently of the degree of the loss of information. The following result describes their pgfs:

**Proposition 1.** *Let $\mathcal{F}$ be a statistical model parameterized by the population mean $\mu$ and described by the set of pgfs $\Phi(t; \mu)$. Then it is closed under binomial subsampling if and only if $\Phi(t; \mu) = g(\mu(t-1))$, for a certain real analytic function $g(x)$ and $\Theta = (0, \mu_M]$, $\mu_M < \infty$ or $\Theta = (0, \infty)$.*

**Proof.** We first prove sufficiency. Let $X$ be a random variable having a pgf of the form $g(\mu_0(t-1))$, where $\mathrm{E}(X) = \mu_0$, $\mu_0 \in \Theta$. Direct calculations from (1) show that the pgf of any $p$-thinning of $X$ is precisely $g(p\mu_0(t-1))$. Consequently, $p\mu_0 \in (0, \mu_M]$ and the model is closed under binomial subsampling.

Turning to necessity, consider a fixed $\mu_0 \in \dot{\Theta}$ (here $\dot{\Theta}$ means the interior of $\Theta$), where $\mu_0 > 0$ because it represents the population mean of a count random variable. Any $0 < \mu < \mu_0$ belongs to $\Theta$ because $\mathcal{F}$ is closed under $p$-thinning and $\mu$ is the expectation of an independent $p$-thinning ($p = \mu/\mu_0$) of a random variable whose pgf is $\Phi(t; \mu_0)$. In terms of their pgfs, by (1) this is equivalent to the identity

$$\Phi(t; \mu) = \Phi\left(1 - \frac{\mu}{\mu_0}(1-t); \mu_0\right). \tag{3}$$

We define $g_{\mu_0}(x) = \Phi(1 + x/\mu_0; \mu_0)$, which is real analytic for $|1 + x/\mu_0| < 1$ or, equivalently, for $-2\mu_0 < x < 0$. Then (3) can be rewritten as

$$\Phi(t; \mu) = g_{\mu_0}(\mu(t-1)),$$

and this is valid for any $0 < \mu < \mu_0$.

Consider now a fixed $\mu_1 \in \dot{\Theta}$, $\mu_1 > \mu_0$. Repeating the preceding argument, defining $g_{\mu_1}(x) = \Phi(1 + x/\mu_1; \mu_1)$, we conclude that $\Phi(t; \mu) = g_{\mu_1}(\mu(t-1))$, for any $0 < \mu < \mu_1$. Consequently $g_{\mu_0}(\mu(t-1)) = g_{\mu_1}(\mu(t-1))$, for $0 < \mu < \mu_0$. Since $g_{\mu_0}(x)$ and $g_{\mu_1}(x)$ are real analytic and they overlap in the open set $-2\mu_0 < x < 0$, we obtain $g_{\mu_0}(x) = g_{\mu_1}(x) = g(x)$ and therefore the function $g(\cdot)$ does not depend on the initial values of $\mu_0$ or $\mu_1$. It follows directly that $\Phi(t; \mu) = g(\mu(t-1))$, for any $\mu \in \dot{\Theta}$, and $\dot{\Theta}$ must be of the form $(0, \mu_M]$ where $\mu_M$ can be $\infty$. In order to conclude that $\Theta = (0, \mu_M]$, it is enough to prove that if $g(\mu(t-1))$ is a pgf, for any $0 < \mu < \mu_M < \infty$, then $g(\mu_M(t-1))$ is also a pgf. This is immediate because $g(x)$ is real analytic in $-2\mu_M < x < 0$. Then it can be expanded in a power series around $-\mu_M$ as

$$g(x) = \sum_{k=0}^{\infty} \frac{g^{(k)}(-\mu_M)}{k!}(x + \mu_M)^k,$$

being convergent for $-\mu_M < x < 0$. Making $x = \mu_M(t-1)$, we obtain

$$g(\mu_M(t-1)) = \sum_{k=0}^{\infty} \frac{g^{(k)}(-\mu_M)}{k!} \mu_M^k t^k,$$

which is real analytic for $-2 < t < 1$. On the other hand, we know that $g(\mu(t-1)) = \sum_{k=0}^{\infty} g^{(k)}(-\mu)\mu^k t^k/k!$ is a pgf for any $0 < \mu < \mu_M$. Consequently, $\sum_{k=0}^{\infty} g^{(k)}(-\mu)\mu^k/k! = 1$ and $g^{(k)}(-\mu)\mu^k/k! \geq 0$ for $k = 0, 1, 2, \ldots$, for any $0 < \mu < \mu_M$. By continuity we conclude that $\sum_{k=0}^{\infty} g^{(k)}(-\mu_M)\mu_M^k/k! = 1$ and $g^{(k)}(-\mu_M)\mu_M^k/k! \geq 0$ for $k = 0, 1, 2, \ldots$, which completes the proof.                                                                       □

The result of this proposition can easily be extended to count models with more than one parameter, such that their vector of parameters has the form $\theta = (\mu, \xi_1, \xi_2, \ldots, \xi_r)$, where any $\xi_i$, $i = 1, \ldots, r$, is invariant under independent $p$-thinning. These $\xi_i$ could, for instance, be the quantities $\mu_{(i)}/\mu^i$ or $\kappa_{(i)}/\mu^i$ specified in Lemma 1 or the quantities in (2). This extension is specified in the following corollary:

**Corollary 1.** *Let $\mathcal{F}$ be a statistical model parameterized by $\theta = (\mu, \xi) \in \Theta \subset \mathbb{R}^r$, where $\mu$ is the population mean and $\xi = (\xi_2, \xi_3, \ldots, \xi_r)$ is invariant under binomial subsampling. Let $\Phi(t; \theta)$ be the set of pgfs that describes the model. Then it is closed under binomial subsampling if and only if $\Phi(t; \theta) = g(\mu(t-1), \xi)$, for certain real analytic functions $g(x, \xi)$. Moreover, fixing $\xi = \xi_0$, the domain of $\mu$ is of the form $(0, \mu_{\xi_0}]$, $\mu_{\xi_0} < \infty$, or $(0, \infty)$.*

*Remark.* Proposition 1 can also be proved by using Theorems 2.4 and 2.6 of Wiuf and Stumpf [18]. Theorem 3.1 of these authors is equivalent to our Corollary 1 when it is assumed that the distributions of the model have finite moments of all orders and these moments determine the distribution. Notice that Corollary 1 is more general.

In order to have a better understanding of the structure of the models which are closed under binomial subsampling, we state the following result that is an extension of Theorem 2.3 of Grandell [5].

**Proposition 2.** *Let $\mathcal{F}$ be a statistical model, closed under binomial subsampling, parameterized by $\theta = (\mu, \xi) \in \Theta \subset \mathbb{R}^r$, where the parameter $\xi$ is also invariant under independent $p$-thinning. Consider a random variable $X$ whose distribution belongs to $\mathcal{F}$; then it satisfies one of the following two assertions:*

  (i) *$X$ is mixed Poisson distributed.*
  (ii) *$X$ is an independent $p$-thinning of a random variable, say $Y$, whose distribution belongs to $\mathcal{F}$ such that $Y$ is a top inverse.*

**Proof.** By Corollary 1, fixing $\xi = \xi_0$, the domain of $\mu$ is of the form $(0, \infty)$ or $(0, \mu_{\xi_0}]$, $\mu_{\xi_0} < \infty$. The two assertions depend on whether the right extreme of the domain of $\mu$ is included or not. Let $X$ be a random variable whose distribution belongs to $\mathcal{F}$ and $\mathrm{E}(X) = \mu_0$:

  (i) Consider that the domain of $\mu$ is of the form $(0, \infty)$. Then $X$ can be obtained by independent $p$-thinnings of random variables whose distributions belong to $\mathcal{F}$, having parameter $\theta' = (\mu_0/p, \xi_0)$ for any $p \in (0, 1]$. Then, by Theorem 2.3 of Grandell [5], $X$ is mixed Poisson distributed.
  (ii) If the domain of $\mu$ is of the form $(0, \mu_{\xi_0}]$, $\mu_0 < \mu_{\xi_0}$, $X$ is an independent $p$-thinning of $Y$ whose distribution has parameter $(\mu_{\xi_0}, \xi_0)$, for $p = \mu_0/\mu_{\xi_0}$. Notice that $Y$ is not an independent $p$-thinning of any random variable in $\mathcal{F}$ except itself. Moreover, this is not an independent $p$-thinning of any random variable because, as the domain $\Theta$ is the widest possible, such a random variable would belong to $\mathcal{F}$. □

What distributions characterize the class of top inverse random variables? The existence of a zero gap is a sufficient but not a necessary condition for a random variable $Y$ to be a top inverse. However, if the distribution of $Y$ has compact support this condition is also necessary (Corollary 3.1 in Thedéen [17]). Example 3 below shows how to construct top inverses without zero-gaps. As a referee has pointed out, further research on characterizations of the class of top inverses would be very interesting. Yannaros [19] has some results related to inverses of thinned renewal processes.

We conclude this section with some examples of models closed under binomial subsampling and show how Propositions 1 and 2 and Corollary 1 can be used.

**Example 1.** This example constitutes an interesting application of Corollary 1. Consider the Hermite distribution introduced by Kemp and Kemp [7, 8]. A convenient parametrization is by means of its population mean and variance. Its pgf has the form

$$\Phi(t; \mu, \sigma^2) = \exp\left\{\frac{\sigma^2 - \mu}{2}(t^2 - 1) + (2\mu - \sigma^2)(t - 1)\right\}.$$

The domain of parameters is $\Theta = \{(\mu, \sigma^2) : \mu > 0, \mu \leq \sigma^2 \leq 2\mu\}$. Changing the parameters to $(\mu, \eta_1)$ by using (2), we obtain

$$\Phi(t; \mu, \eta_1) = \exp\left\{\frac{\mu^2 \eta_1}{2}(t - 1)^2 + \mu(t - 1)\right\},$$

where the domain of parameters is now $\Theta' = \{(\mu, \eta_1) : \eta_1 \geq 0, 0 < \mu \leq 1/\eta_1\}$. Notice that $\Phi(t; \mu, \eta_1) = g(\mu(t-1), \eta_1)$, where $g(x; \eta_1) = \exp(\eta_1 x^2/2 + x)$. Consequently the Hermite family is closed under binomial subsampling. It is known that the Hermite distribution is not a mixed Poisson model (Kemp and Kemp [7]). It can easily be shown that any Hermite random variable with fixed $\mu_0$ and $\eta_{1_0}$ can be obtained by independent $p$-thinning from another random variable with pgf

$$\Phi^*(t) = \exp\left\{\frac{1}{2\eta_1}(t^2 - 1)\right\}$$

with $p = \mu\eta_1$. This pgf corresponds to a random variable of the form $2X$, where $X$ is Poisson distributed, and this is a top inverse random variable because $P(2X = k) = 0$ for $k = 1, 3, 5, \ldots$. It is in concordance with Proposition 2.

**Example 2.** Consider now the negative binomial family of distributions parameterized by its population mean and variance. Its pgf has the form

$$\Phi(t; \mu, \sigma^2) = \exp\left(-\frac{\mu^2 \log(1 - (\sigma^2/\mu - 1)(t - 1))}{\sigma^2 - \mu}\right).$$

The domain of parameters is $\Theta = \{(\mu, \sigma^2) : \mu > 0, \mu \leq \sigma^2\}$. Changing the parameters to $(\mu, \eta_1)$ as in Example 1, we obtain

$$\Phi(t; \mu, \eta_1) = (1 - \mu\eta_1(t - 1))^{-1/\eta_1},$$

where the domain of parameters is now $\Theta' = \{(\mu, \eta_1) : \mu > 0, \eta_1 \geq 0\}$. Notice that $\Phi(t; \mu, \eta_1) = g(\mu(t - 1), \eta_1)$, where $g(x; \eta_1) = (1 - \eta_1 x)^{-1/\eta_1}$. Corollary 1 again shows that the negative binomial family of distributions is closed under binomial subsampling.

This family is a mixed Poisson model. Notice that any negative binomial random variable with parameters $(\mu_0, \eta_{1_0})$ is an independent $p$-thinning of a negative binomial random variable with parameters $(\mu_0/p, \eta_{1_0})$, for any $p \in (0, 1]$.

**Example 3.** Consider a count random variable $X$ such that $p_k = P(X = k) > 0$ for $k = 0, 1, 2, \ldots$ and $p_{2k+1}/p_{2k}$ (or $p_{2k}/p_{2k-1}$) tends to infinity when $k \to \infty$. Our aim is to prove that $X$ is a top inverse. Suppose that there exist a random variable $Y$ such that $X$ can be obtained as an independent $p$-thinning of $Y$; we will arrive at a contradiction. According to (1), the pgf of $Y$ is $\Phi_Y(t) = \Phi_X(1 + (t-1)/p)$. Consequently, the probabilities $p_k^* = P(Y = k)$, $k = 0, 1, 2, \ldots$, satisfy

$$p_k^* = k! \Phi_Y^{(k)}(0) = \frac{k!}{p^k} \Phi_X^{(k)}\left(1 - \frac{1}{p}\right) = \frac{1}{p^k} \sum_{i=0}^{\infty} \binom{i+k}{k} \left(1 - \frac{1}{p}\right)^i p_{k+i}.$$

Adding each pair of consecutive terms, we arrive at the expression

$$p_k^* = (1 + \varepsilon)^k \sum_{i=0}^{\infty} \binom{2i+k}{k} \varepsilon^{2i} p_{2i+k} \left(1 - \varepsilon \frac{2i+k+1}{2i+1} \frac{p_{2i+k+1}}{p_{2i+k}}\right),$$

where $\varepsilon = -1 + 1/p$. Therefore, for any $\varepsilon > 0$, there exists a $r$ such that, for all $i \geq 0$, $p_{2i+r+1}/p_{2i+r} > 1/\varepsilon$. Consequently $p_r^* < 0$, which is a contradiction.

Now many examples of top inverses without zero gaps can be constructed by alternating two series of positive terms. For instance, this is the case for the top inverse random variable $X$ defined by the probability function $p_{2k} = c(1/3)^k$, $p_{2k+1} = c(1/2)^k$, $k = 0, 1, 2, \ldots$, where $c = 2/7$ is the normalizing constant. Its pgf is $\Phi_X(t) = 6/(21 - 7t^2) + 4t/(14 - 7t^2)$. According to (1), the statistical model constructed by binomial subsampling of $X$ is described by the pgfs $\Phi(t; p) = 6/(21 - 7(1 - p(1 - t))^2) + 4(1 - p(1 - t))/(14 - 7(1 - p(1 - t))^2)$, where $p \in (0, 1]$. It is immediately evident that the population mean is $\mu = 15p/7$ and, consequently, the model parametrized by $\mu$ is described by the pgfs $\Phi(t; \mu) = g(\mu(t - 1))$, where $g(x) = (900 + 420x)/(3150 - 7(15 + 7x)^2) + 1350/(4725 - 7(15 + 7x)^2)$. Proposition 1 shows that this model is closed under binomial subsampling. Notice that the domain of parameters is $\mu \in (0, 15/7)$ and cannot be extended. This is again in concordance with Proposition 2.

## 3. Characterization

In many applications it is reasonable to use count models which have the property of conservation under addition of independent effects, that is, their distributions are closed under convolutions.

**Definition 3.** *Let $\mathcal{F}$ be a statistical model. It is said to be closed under addition if, for any two independent random variables $X$, $Y$ with distributions belonging to $\mathcal{F}$, the distribution of $X + Y$ also belongs to $\mathcal{F}$.*

For instance, this property is satisfied by the Poisson model or by the Hermite family of distributions (see Puig [13]). A quite natural parameterization of a model closed under addition is by means of its cumulants. This is due to the fact that the $i$th cumulant of the sum of two independent random variables is the sum of the corresponding $i$th cumulants. It also occurs with factorial cumulants, which therefore provide another reasonable parameterization of a model closed under addition which is convenient for our purposes. The following proposition can be proved immediately using a methodology similar to that of Teicher [16] concerning the structure of the characteristic function of a model closed under addition.

**Proposition 3.** *Let $\mathcal{F}$ be a count statistical model parameterized by its $r$ first factorial cumulants $\theta = (\kappa_{(1)}, \kappa_{(2)}, \ldots, \kappa_{(r)})$, with a pgf continuous in $\theta \in \Theta \subset \mathbb{R}^r$. The model is closed under addition if and only if its pgf can be expressed as*

$$\Phi(t; \theta) = \exp\left\{\sum_{i=1}^{r} \alpha_i(t)\kappa_{(i)}\right\}, \tag{4}$$

*where $\alpha_i(t)$ are functions not depending on the parameters.*

**Proof.** It can immediately be seen that the form of the pgf (4) implies that the model is closed under addition.

Suppose now that the model is closed under addition. This means that, for any pair of values $\theta_1, \theta_2 \in \Theta$, we have

$$\Phi(t; \theta_1 + \theta_2) = \Phi(t; \theta_1)\Phi(t; \theta_2).$$

Taking logarithms and fixing $t$, we obtain

$$f(\theta_1 + \theta_2) = f(\theta_1) + f(\theta_2), \tag{5}$$

where $f(\theta) = \log(\Phi(t; \theta))$. But (5) is the multidimensional Cauchy functional equation. It is known that its solution has the form, $f(\theta) = \sum_{i=1}^{r} \alpha_i \theta_i$, where the $\alpha_i$ do not depend on $\theta$. However, they can depend on $t$, and this leads to the solution $\log(\Phi(t; \theta)) = \sum_{i=1}^{r} \alpha_i(t)\theta_i$, thus concluding the proof. $\qquad \square$

Proposition 1 and Corollary 1, together with Proposition 3, allow us to find all the count statistical models satisfying mild conditions, such that they are closed under addition and binomial subsampling. Surprisingly there do not exist many types of models satisfying both properties.

**Theorem 1.** *Let $\mathcal{F}$ be a count statistical model parameterized by its $r$ first factorial cumulants $\theta = (\kappa_{(1)}, \kappa_{(2)}, \ldots, \kappa_{(r)})$. Assume that the model has a pgf continuous in $\theta \in \Theta \subset \mathbb{R}^r$. The model is closed under addition and under binomial subsampling if and only if it has a pgf of the form*

$$\Phi(t; \theta) = \exp\left\{\sum_{i=1}^{r} \frac{\kappa_{(i)}}{i!}(t-1)^i\right\}. \tag{6}$$

**Proof.** We first prove sufficiency. The form of the pgf (6) implies, by Proposition 3, that the model is closed under addition. Notice that this pgf can also be reparameterized as

$$\Phi(t; \mu, \xi) = \exp\left\{ \mu(t-1) + \sum_{i=2}^{r} \xi_i \mu^i \frac{(t-1)^i}{i!} \right\},$$

where $\xi = (\xi_2, \ldots, \xi_r)$ and $\xi_i = \kappa_{(i)}/\mu^i$. Then the pgf can be written as $\Phi(t; \mu, \xi) = g(\mu(t-1); \xi)$, where $g(x; \xi) = \exp(x + \sum_{i=2}^{r} \xi_i x^i / i!)$. Using Corollary 1, we conclude that the model is closed under binomial subsampling.

Turning to necessity, suppose that the model is closed under addition and under binomial subsampling. Again using the parameterization $(\mu, \xi)$, where $\xi_i = \kappa_{(i)}/\mu^i$, by Proposition 3 and Corollary 1 we obtain the identity

$$\exp\left\{ \alpha_1(t)\mu + \sum_{i=2}^{r} \alpha_i(t)\mu^i \xi_i \right\} = g(\mu(t-1); \xi).$$

Fixing $\xi = \xi_0$ and $t = t_0$ and taking $x = \mu(t_0 - 1)$, we obtain

$$\log(g(x; \xi_0)) = \frac{\alpha_1(t_0)}{t_0 - 1} x + \sum_{i=2}^{r} \frac{\alpha_i(t_0)}{(t_0 - 1)^i} \xi_{0_i} x^i,$$

that is, a polynomial of degree $r$ in $x$. Since it must not depend on the choice of $t_0$, we obtain that the functions $\alpha_i(t)$ in (4) have the form $\alpha_i(t) = c_i(t-1)^i$, for $i = 1, \ldots, r$, where the $c_i$ are constants. But these constants are $c_i = 1/i!$, because the Taylor series coefficients of $\log(\Phi(t; \theta))$ at $t = 1$ must be exactly $\kappa_{(i)}/i!$, $i = 1, \ldots, r$. This concludes the proof. □

The preceding theorem shows that, under mild conditions, the only statistical model closed under addition and under binomial subsampling, such that it can be parameterized by its population mean ($r = 1$), is the Poisson model.

Kemp and Kemp [8] present an elementary proof showing that the only count distribution such that its pgf has the form (6), having a second-degree polynomial in the exponent, is the Hermite distribution. Consequently, the two-parameter case ($r = 2$) described by Theorem 1 corresponds to the Hermite family of distributions.

Notice that the polynomial of degree $r$ that appears in the exponent of (6) can always be expressed as $\sum_{i=1}^{r} a_i(t^i - 1)$. This leads to an alternative parameterization that allows us to represent the pgf of the count models characterized by Theorem 1 in the form

$$\Phi(t; a) = \exp\left\{ \sum_{i=1}^{r} a_i(t^i - 1) \right\}. \tag{7}$$

This representation is interesting because, if the $a_i$ are non-negative, any count random variable with a pgf as in (7) can be understood as a linear combination of independent

Poisson random variables (Feller [2], pages 291–292). Specifically, it can easily be shown that any count random variable $X = \sum_{i=1}^{r} i X_i$, where $X_i$ are independent Poissons with mean $a_i$, has a pgf of the form (7). However, it is possible for $\Phi(t; a)$ to be a pgf even when some of the $a_i$ in (7) are negative, but in this case these pgfs are not infinitely divisible (see Lévy [9]; Lukacs [10], page 251; Milne and Westcott [11]).

Milne and Westcott [11] designate these distributions as *generalized Hermite* or *rth-order univariate Hermite*. We prefer the second name because of Gupta and Jain's [6] use of the name 'generalized Hermite distribution', related to the specific case in (7) where only $a_1$ and $a_r$ are non-zero.

## 3.1. The $r$th-order univariate Hermite distribution

Let $X$ be a count random variable with a pgf as in (7). Consider its probabilities, $p_k = P(X = k)$, $k = 0, 1, \ldots$. Notice that $\Phi(t; a) = \sum_{k=0}^{\infty} p_k t^k = \exp\{\sum_{i=1}^{r} a_i(t^i - 1)\}$. Then, taking logarithms and differentiating with respect to $t$, we obtain the identity

$$p_1 + 2p_2 t + 3p_3 t^2 + \cdots = (p_0 + p_1 t + p_2 t^2 + \cdots) \sum_{i=1}^{r} i a_i t^{i-1}.$$

Equating the terms with the same degree in $t$ in both parts of the identity, we obtain the following lemma:

**Lemma 2.** *Let $X$ be a count random variable with a pgf as in (7). Then its probabilities can be calculated from the recurrence relation $p_k = (\sum_{i=1}^{r} i p_{k-i} a_i)/k$, $k = 1, 2, \ldots$, with $p_0 = \exp\{-\sum_{i=1}^{r} a_i\}$ and $p_{-1} = p_{-2} = \cdots = p_{1-r} = 0$.*

For the case $r = 2$ (Hermite distribution), this particular recurrence relation was found by Kemp and Kemp [7]. Notice that Lemma 2 is a version of the Panjer recursion formula (Panjer [12]).

The parameterization in terms of the $a_i$ is particularly worthwhile in order to calculate the probabilities. Moreover, it can be easily shown that the cumulants can be directly calculated from the $a_i$ by means of the expression, $\kappa_s = \sum_{i=1}^{r} i^s a_i$. When all the $a_i$ are non-zero, it is known (see Remark 1 in Milne and Westcott [11]) that if at least $a_1$, $a_{r-1}$ and $a_r$ are non-zero, then $a_1 > 0$, $a_{r-1} > 0$ and $a_r > 0$ are necessary conditions for (7) to be a pgf. Consequently, the full domain of parameters for $r = 2$ and $r = 3$ does not include negative values for any $a_i$ but, for $r \geq 4$, some of these coefficients could be negative. The domain of the parameters that will be considered subsequently for practical applications is $\Theta = \{a = (a_1, \ldots, a_r) : a_i \geq 0\}$. This ensures that the distributions are infinitely divisible. Moreover, when all the $a_i$ are non-negative, it can easily be checked that the distribution is overdispersed, that is, the variance of the count variable exceeds its mean. Consequently, these models can be used in practice to analyse data sets when we encounter the phenomenon of overdispersion.

In order to choose the appropriate order $r$ of the model, a hierarchical approach based on the likelihood ratio test can be used. In our context this means testing the hypothesis $H_0 : a_r = 0$. It is important to take into account that, under the null hypothesis, the likelihood ratio test statistic does not have an asymptotic $\chi_1^2$ distribution as usually happens, because $a_r = 0$ belongs to the boundary of the domain of parameters. It can be established that in this situation the asymptotic distribution of the likelihood ratio test statistic is a 50:50 mixture of the constant zero and the $\chi_1^2$ distribution (Self and Liang [15]; Geyer [3]). The $\alpha$ upper tail percentage points for this mixture are the same as the $2\alpha$ upper tail percentage points for a $\chi_1^2$ distribution.

# Acknowledgements

The authors thank the editor, associate editor and a referee for their valuable comments and suggestions. This research was partially supported by grant MTM2006-01477 from the Ministry of Education of Spain.